\documentclass{amsart}
\usepackage{amsfonts}
\usepackage{amsmath,amscd}
\usepackage{amssymb}
\usepackage{amsthm}
\usepackage{newlfont}
\usepackage{graphicx}
\usepackage[all]{xy}
\usepackage[usenames]{color}
\usepackage{comment}

\long\def\alert#1{\parindent2em\smallskip\hbox to\hsize
{\hskip\parindent\vrule
\vbox{\advance\hsize-2\parindent\hrule\smallskip\parindent.4\parindent
\narrower\noindent#1\smallskip\hrule}\vrule\hfill}\smallskip\parindent0pt}
 \newtheorem{thm}{Theorem}[section]

 \newtheorem{lem}[thm]{Lemma}
 \newtheorem{prop}[thm]{Proposition}
\theoremstyle{definition}
 \newtheorem{defn}[thm]{Definition}
\theoremstyle{remark}

 \theoremstyle{problem}
 
 \numberwithin{equation}{section}
\newtheorem*{Theorem A}{\textbf{Theorem A}}
\newtheorem*{Corollary E}{\textbf{Corollary E}}
\newtheorem*{Corollary F}{\textbf{Corollary F}}
\newtheorem*{Theorem B}{\textbf{Theorem B}}
\newtheorem*{Theorem C}{\textbf{Theorem C}}
\newtheorem*{Corollary D}{\textbf{Corollary D}}
\newtheorem*{p a}{\textbf{Proof of  Theorem A}}
\newtheorem*{p b}{\textbf{Proof of Theorem B}}
\newtheorem*{p c}{\textbf{Proof of  Theorem C}}
\newtheorem*{p d}{\textbf{Proof of  Corollary D}}
\newtheorem*{p e}{\textbf{Proof of  Corollary E}}
\begin{document}

%
\title[Classification of  some nilpotent Lie algebras ]{Classification of nilpotent Lie algebras of nilpotency class 3 having a derived subalgebra of dimension three}
\author{Saboura Yousefi, Azam Kaheni and Farangis Johari}


\begin{abstract}
In this paper, we investigate nilpotent Lie algebras  $ L $  of nilpotency class $3 $ and provide a complete classification of those satisfying  $ \dim L^2 = 3 $ and  $Z(L) = L^3 \cong A(2). $  Furthermore, we explicitly characterize the structure of such Lie algebras in the case when $ \dim L = 7. $

\end{abstract}
\maketitle

\maketitle
\section{Introduction}
The classification of Lie algebras is a fundamental problem that has long attracted substantial attention, and one fundamental invariant in this context is the dimension of the Lie algebra.
 In 1977, Skjelbred and Sund introduced a systematic method for constructing all nilpotent Lie algebras of dimension \( n \) from the classification of those of smaller dimension together with their automorphism groups. Since then, numerous efforts have advanced this line of research, resulting in several classification lists. In dimension \(7\), multiple independent classifications have appeared. Saullina (1964) provided the first list over \( \mathbb{C} \), which was subsequently extended by Romdhani to cover both \( \mathbb{R} \) and \( \mathbb{C} \) \cite{24,25}. Later, Seeley~\cite{31} and Ancochea--Goze~\cite{2} obtained further classifications over \( \mathbb{C} \). These lists rely on different invariants and consequently exhibit certain discrepancies.

In 1989, Carles introduced the weight system as a new invariant and compared the classifications of Saullina, Romdhani, and Seeley, identifying several omissions and inaccuracies. Building on Carles’ observations and earlier results, Seeley produced a revised and consolidated classification over \( \mathbb{C} \) \cite{33}. Gong~\cite{gon} later employed the Skjelbred--Sund construction to obtain all \(7\)-dimensional nilpotent Lie algebras, advocating this approach for its uniformity and completeness. Additional partial classifications have been developed for particular families of nilpotent Lie algebras, including nilpotent algebras of maximal rank (Favre~\cite{10}), two-step algebras (Scheuneman~\cite{30}, Gauger~\cite{11}, Revoy~\cite{23}), and filiform algebras (Ancochea--Goze~\cite{1}).
Among the structural invariants that play a significant role in the analysis of Lie algebras, the Schur multiplier is of particular importance. For nilpotent Lie algebras, several equivalent formulations of this invariant are known. In particular, Batten showed that the Schur multiplier of a nilpotent Lie algebra \(L\) admits the description
\begin{align}\label{eq1}
    M(L)\cong (R \cap \gamma_{2}(F)) / [R,F],
\end{align}

where $L\cong F/R $  and \(\gamma_{2}(F)\) denotes the derived subalgebra of  free Lie algebra \(F\) (\cite[Theorems~1.12 and~3.6]{batten}).

In 1996, Batten \emph{et al.} introduced the invariant \( t(L) \) and established a classification framework based on the dimension of a Lie algebra and its Schur multiplier \cite{bms,H,HS}. 
The upper bound in \eqref{eq1} naturally leads to the definition of the numerical invariant, \( t(L) \), given by
\begin{align}\label{eq2}
    t(L) = \binom{n}{2} - \dim M(L),
\end{align}
where \(n\) denotes the dimension of the Lie algebra \(L\).
In 1998, Ming-peng classified all finite-dimensional Lie algebras up to dimension seven, determining their explicit structures \cite{gon}. For dimensions greater than \(7\), only certain families with particular structural properties have been completely determined. For instance, Lindsey \cite{B} investigated Lie algebras of maximal class and classified those satisfying $t(L) < 17.$

In 2011, Niroomand \emph{et al.} refined this framework by introducing the invariant \( s(L) \), which incorporates the dimension of the derived subalgebra into \( t(L) \). Their work expanded the classification obtained via \( t(L) \), yielding precise structural descriptions of additional Lie algebras \cite{N,NS2022,SAN,NS2020}. More recently, Johari \emph{et al.} examined nilpotent Lie algebras by simultaneously considering their nilpotency class and the dimension of \( L^{2} \), determining several remaining structures \cite{FAP,AFP}.

An upper bound for the dimension of the derived subalgebra, expressed in terms of the dimension of the quotient by the center, was established by Moneyhun in \cite[Lemma~1]{1991}. In particular, for a Lie algebra \(L\) satisfying 
\[
\dim\!\left( \frac{L}{Z(L)} \right) = n,
\]
the following inequality is obtained:
\begin{align}\label{eq6}
    \dim \gamma_{2}(L) \leq \binom{n}{2}.
\end{align}

Motivated by the bound in~\eqref{eq6}, we introduce a new invariant for nilpotent Lie algebras, denoted by \( t(L^{2}) \), formulated in terms of the dimension of the derived subalgebra. This invariant measures the discrepancy between the maximal achievable dimension of the derived subalgebra determined by the codimension of the center and its actual dimension. Thus, it serves as a complementary tool to existing invariants such as \( t(L) \) and \( s(L) \).

\begin{defn}
Let \( L \) be a nilpotent Lie algebra with \( \dim\!\left( \frac{L}{Z(L)} \right) = n \). We define
\begin{align}
    t(L^{2}) = \frac{1}{2} n(n - 1) - \dim \gamma_{2}(L).
\end{align}  
\end{defn}
Furthermore, Batten {\it et al.} ~\cite{bms} provided a classification of the central factor of Lie algebras $L$ according to the invariant $t(L^2),$ which takes values in $\{0,1\}.$
We first established a connection between the invariant \( t(L^{2}) \), the dimension of a central factor of \(L\), the dimension of \(L\) itself, and the dimension of its quotient by the derived subalgebra. In previous work \cite{ukj}, we provided a complete classification of finite-dimensional nilpotent Lie algebras whose central factor has dimension \(n\) in the cases \( t(L^{2}) = 0 \), \(1\), and \( \tfrac{1}{2}n(n-1)-1 \) (Theorems~A, B, and~C).

When \( t(L^{2}) = \tfrac{1}{2}n(n-1)-3 \) and \( \dim Z(L)=1 \), the Lie algebra must have nilpotency class $3.$
The finite-dimensional nilpotent Lie algebras satisfying these conditions were classified in the following result:

\begin{thm}\cite[Main Theorem]{AFP}\label{Z(L)=1}
Let \( L \) be a finite-dimensional nilpotent Lie algebra of class \(3\) with \(\dim L^{2} = 3\) and \( \dim Z(L) = 1\). Then \(L\) is isomorphic to one of the following non-isomorphic Lie algebras:
\[
L_{6,19(\epsilon)},\; L_{6,20},\quad \text{and} \quad S_{1}, S_{2},\dots, S_{32}.
\]
\end{thm}

We are motivating to give  complete classification of finite-dimensional nilpotent Lie algebras whose central factor has dimension \( 5 \) and nilpotency class is 3,  $t(L^2) = \frac{1}{2}n(n-1)-3$ when \(L^{3}=Z(L)\cong A(2)\). It follows immediately that the dimension of any such Lie algebra is $7.$ 

Throughout this paper, all Lie algebras are assumed to be finite-dimensional over a field \( \mathbb{F} \) of characteristic not equal to \(2\), and the Lie bracket is denoted by \([\, , ]\). The purpose of this work is to classify all nilpotent Lie algebras of class \(3\) for which the derived subalgebra has dimension three and the third term of the lower central series coincides with the center, that is, \(L^{3}=Z(L)\cong A(2)\). In particular, we obtain a complete description of their structure in arbitrary dimension and, as a concrete application, provide the explicit classification of all such Lie algebras in dimension \(7\).

     \section{Preliminaries}
    
    In this section, we recall several essential and classical results that will play a key role in the subsequent analysis. \\
 A Lie algebra $L$ is called a semidirect sum of an ideal $I$ by a subalgebra $J$ if $L=I+J$ and $I \cap J=0.$ In this case, we denote it by $L=I\rtimes J.$\\
 An abelian Lie algebra of dimension $n$ is denoted by $A(n).$ A Lie algebra $H$ is called a generalized Heisenberg Lie algebra of rank $n\geqslant 1$, if $H^{2}=Z(H)$ and $dim~H^{2}=n$. In particular, when $n=1$, $H$ is the classical Heisenberg Lie algebra. Furthermore, a Lie algebra $L$ is called a stem Lie algebra if  if its center $Z(L)$ is contained in its derived subalgebra $L^{2},$  that is, $Z(L)\subseteq L^{2}.$
It is well known \cite [Lemma 2.3]{Johari2022}, that for any non-abelian Lie algebra $L$ with finite $dim~\dfrac{L}{Z(L)},$ there exists a decomposition
\[
L\cong K\oplus A,~L^{2}\cong K^{2}\quad\text{and}\quad Z(K)\cong L^{2}\cap Z(L)
\]
where $A $ is an abelian ideal.These observations motivate the following result.
\begin{lem}\label{l1}
Let $L$ be a nonabelian Lie algebra such that $dim~\dfrac{L}{Z(L)}$ is finite. Then $L\cong S\oplus A$, where $S$ is a stem Lie algebra and $A$ is an abelian Lie algebra.
\end{lem}
The concept of the Schur multiplier plays a fundamental role in precisely characterizing the structure of Lie algebras.
Let $L$ be a nilpotent Lie algebra of dimension $n$ in a field $F,$ and consider the central extension  
\[
0 \longrightarrow M \longrightarrow C \longrightarrow L \longrightarrow 0.
\]
  If $M$ attains the maximal possible dimension, then $M\cong M(L).$  Consequently, we have the inequality $dim~Z(L)\cap L^{2} \leqslant dim~M(L).$

The nilpotency class of a Lie algebra \( L \), denoted by \( \text{cl}(L) \), is defined as the smallest positive integer \( m \) such that \( L^{m+1} = 0 \). If no such $m$ exists, \( L \) is said to be non-nilpotent.
\begin{prop}\label{class}\cite [Lemma 2.5]{FAP}
Let $L$ be an $n$-dimensional Lie algebra of nilpotency class \(c.\) Then\\
\[
 L^{c-i}\nsubseteq Z_{i}(L)~~\text{for ~all}~~0\:<i\le\:c-1.
 \]
\end{prop}
The following result describes the relationship between the terms of the lower central series of a Lie algebra and their corresponding subalgebras.
\begin{lem}\cite[Lemma1]{z}\label{low}
Let $L$ be a nilpotent Lie algebra, and let $K$ be a subalgebra of $L$ such that $L^{2}=K^{2}+L^{3}.$ Then, for $i\geqslant2$, 
\[
L^{i}=K^{i}+L^{i+1},\quad \text{and also}\quad L^{i}=K^{i}.
\]
Moreover, $K$ is an ideal of $L.$
\end{lem}
A Lie algebra is called capable if it can be realized as a central quotient of another Lie algebra.
The precise structure of any capable nilpotent Lie algebra whose derived subalgebra has dimension one is given as follows:
\begin{thm}\cite[Theorm 3.6]{Niroo2013}\label{Theorem 3.6}
    Let $L$ be a nilpotent Lie algebra of dimension $n$ such that $dim~L^2=1.$ Then $L\cong H(m)\oplus A(n-2m-1)$ for some $m$ and $L$ is capable if and only if $m=1,$ that is $L\cong H(1)\oplus A(n-3).$
\end{thm}

In the following, we present the exact structure of an $n$-dimensional nilpotent stem Lie algebra of nilpotency class $3$ whose derived subalgebra has dimension $3.$
\begin{lem}\cite[Lemma 2.1]{AFP}\label{Lemma 2.1}
Let $T$ be an n-dimensional nilpotent stem Lie algebra of class $3$ with derived subalgebra of dimension $3.$ Then:
\begin{enumerate}
\item[$(i)$]
     If $T^3=Z(T)\cong A(1),$ Then $ \frac{T}{Z(T)}\cong L_{5,8}\oplus A(n-6), L_1 \oplus A(n-8),$ or $L_{6,22(\epsilon)}\oplus A(n-7) $ for all $\epsilon \in \frac{\mathbb{F}}{\sim^*}.$  
\item[$(ii)$]  If  $T^3\cong A(1)$ and $ Z(T)\cong A(2), $  then $\frac{T}{Z(T)}\cong H(1)\oplus A(n-5).$
   \item[$(iii)$] If $T^3= Z(T)\cong A(2), $ then $ \frac{T}{Z(T)}\cong H(1)\oplus A(n-5).$
\end{enumerate}
\end{lem}
Recall that, the Lie algebra $L $ is a central product of  two ideals $A$ and $B,$ if $L = A + B,$
 such that  $[A, B] = 0$  and $A \cap B \subseteq Z(L). $ We denote the central product of two Lie algebras $A$ and $B $ by \( A \dot{+} B \). \par 



\subsection{Main Results} \label{seca}
We now concentrate on nilpotent Lie algebras of class $ 3$ and provide a complete classification of those with  $\dim L^2 = 3. $ \par By Lemma \ref{l1} it suffices to focus on stem Lie algebras satisfying $\dim L^2 = 3.$ Furthermore, Proposition \ref{class} ensures that $L^2 \nsubseteq Z(L),$  which restricts our attention to Lie algebras where $ \dim Z(L) = 1$ or $ \dim Z(L) = 2.$
The case $ \dim Z(L) = 1,$ has been completely settled and is given in Theorem \ref{Z(L)=1}.
The following theorem describes the structure of all finite-dimensional nilpotent Lie algebras \( L \) of nilpotency class \( 3 \) with a derived subalgebra of dimension \( 3 \), under the condition that \( \dim Z(L) = \dim \gamma_3(L) = 2 \). 
\begin{Theorem A}
{\em
Let \( L \) be an \( n \)-dimensional nilpotent stem Lie algebra over a field \( \mathbb{F} \) with \( \mathrm{char}(\mathbb{F}) \neq 2 \), satisfying  $cl(L)=3,$  \( \dim Z(L) = \dim \gamma_3(L) = 2 \) and \( \dim \gamma_2(L) = 3 \). Suppose \( I \), \( A \), \( T \), and \( K \) are subalgebras of \( L \) such that  \( I \cong L_{5,9}=\langle x_i\mid [x_1, x_2] = x_3,
    [x_1, x_3] = x_4, [x_2, x_3] = x_5\rangle \), \( T \cong H(m) \), \( A \) is abelian, and \( K \) is a generalized Heisenberg Lie algebra of rank two.
Then \( L \) is isomorphic to one of the following structures:
\begin{enumerate}
\item[$(i)$]  \( L = I \rtimes A \), where \( \dim A = n-5 \).
\item[$(ii)$] \( L = I \dot{+} T \), where \( n = 2m + 5 \) for some \( m \geq 1 \).
\item[$(iii)$] \( L = (I \rtimes A) \dot{+} T \), where \( A \cong A(n-2m-5) \) and $ 0 \neq [I, A] \subseteq Z(L).$ 
\item[$(iv)$] \( L = I + T \), where $ [I, T]=Z(L)$ and \( n = 2m + 5 \) for some \( m \geq 1 \).
\item[$(v)$] \( L = (I \rtimes A) + T \), where $0\neq [I,A] \subseteq [I, T]=Z(L)$, \( A \cong A(n-2m-5) \). 
\item[$(vi)$] \( L = I + T \), where $ [I, T]\cong A(1)$ and \( n = 2m + 5 \) for some \( m \geq 1 \). If $[I, T]=T^2,$ then $T$ is ideal.
\item[$(vii)$] \( L = (I \rtimes A) + T=(I+T)\rtimes A \), where $0\neq [I,T]\cong A(1) \subseteq Z(L)$ and $0\neq [I,A]\subseteq Z(L)$ , \( A \cong A(n-2m-5) \).
\item[$(viii)$] Either \( L = I \dot{+} K \), where \( \dim K = n-3 \geq 5 \) or $L = I + K$, where $0\neq [I ,K]\subseteq Z(L)=K^2.$ Moreover, $K$ is ideal.
\item[$(ix)$]Either \( L = (I \rtimes A) \dot{+} K \) where $ 0 \neq [I, A] \subseteq Z(L),$  \( A \cong A(r) \) and \( \dim K = n - r - 3 \) for some \( r \geq 1 \) or \( L = (I \rtimes A) + K \) where $ 0 \neq [I, A] \subseteq Z(L),$ $0\neq [I ,K]\subseteq Z(L)=K^2,$ \( A \cong A(r) \) and \( \dim K = n - r - 3 \) for some \( r \geq 1 \).
\end{enumerate}
Whenever the corresponding subalgebras appear, the following conditions are satisfied:
$
\gamma_2(I) = \gamma_2(L),$ $  \gamma_2(K)=Z(L) = Z(I)=\gamma_3(I) = \gamma_3(L), $   and $  \gamma_2(T) \subsetneq Z(L).
$}
\end{Theorem A}

We now provide a full classification of all non-isomorphic $7$-dimensional stem nilpotent Lie algebras over any field with characteristic different from two, under the given conditions.

\begin{Theorem B}
{\em
Let $L$ be a $7$-dimensional stem nilpotent Lie algebra  over a field \( \mathbb{F} \) with \( \mathrm{char}(\mathbb{F}) \neq 2 \), satisfying  $cl(L)=3,$  $\dim \gamma_2(L)=3,$ and $Z(L)=\gamma_3(L)\cong A(2)$. Then $L$ is isomorphic to one of the following non-isomorphic Lie algebras:
\[
L_{1}=\langle x_{i}\mid [x_{1},x_{2}]=x_{3},~
[x_{1},x_{3}]=[x_{1},x_{4}]=x_{6},~[x_{2},x_{3}]=[x_{1},x_{5}]=x_{7},1\leq i\leq 7\rangle.    
\]
\[
L_{2}=\langle x_{i}\mid [x_{1},x_{2}]=x_{3},~
[x_{1},x_{3}]=[x_{2},x_{4}]=x_{6},~[x_{2},x_{3}]=[x_{1},x_{5}]=x_{7},~1\leq i\leq 7\rangle.   
\]
$
L_{3}=\langle  x_{i}\mid[x_{1},x_{2}]=x_{3},~
[x_{1},x_{3}]=x_{6},~[x_{2},x_{3}]=[x_{4},x_{5}]=x_{7},~1\leq i\leq 7\rangle.
$
\[
L_{4}=\langle  x_{i} \mid[x_{1},x_{2}]=x_{3},~
[x_{1},x_{3}]=[x_{2},x_{4}]=x_{6},~[x_{2},x_{3}]=[x_{4},x_{5}]=x_{7},~1\leq i\leq 7\rangle.
\]
}
\end{Theorem B}
 
\section{Proofs of main results}
This section is devoted to establishing the proofs of the main results stated in Subsection~\ref{seca}. We are now ready to present the proof of Theorem A, which follows an approach similar to that used in \cite[Proposition 2.4]{AFP} and \cite[Proposition 4.3 and Theorem 4.6]{NP1}.

\begin{p a}
{\em
Let \( n \geq 6 \). Using Lemma \ref{Lemma 2.1}, we have \(\dfrac{L}{Z(L)} \cong H(1) \oplus A(n-5)\). Therefore, there exist two ideals \(\dfrac{I_{1}}{Z(L)}\) and \(\dfrac{I_{2}}{Z(L)}\) of \(\dfrac{L}{Z(L)}\) such that 
\[
\dfrac{I_{1}}{Z(L)} \cong H(1) \quad \text{and} \quad \dfrac{I_{2}}{Z(L)} \cong A(n-5).
\]
Since \(L^2 / Z(L) = (I_1^2 + Z(L)) / Z(L)\) and \(Z(L) = L^3\), by  Lemma \ref{low}, we obtain \(L^2 = I_1^2\) and \(L^3 = I_1^3\), so \(c(L) = c(I_1) = 3\).

We claim that \(Z(L) = Z(I_1)\).  
Since \(Z(I_1)/Z(L) \subseteq Z\left(\dfrac{I_{1}}{Z(L)}\right) = \dfrac{I_{1}^2}{Z(L)}\) and \(Z(I_1) \neq I_1^2\), we have
\[
A(2) \cong Z(L) \subseteq Z(I_1) \subseteq I_1^2 \cong A(3),
\]
and thus \(Z(L) = Z(I_1)\).  It is noteworthy that, according to the classification of $5$-dimensional Lie algebras in \cite{graf}, the unique $5$-dimensional nilpotent stem Lie algebra $L$ of class $3$ with $\dim L^2=3$ and $\dim Z(L)=\dim \gamma_3(L)=2$ is $L_{5,9}.$ Consequently, we have \(I_1 \cong L_{5,9}\).

Since \(I_1\) and \(I_2\) are ideals of \(L\), their commutator satisfies \([I_1, I_2] \subseteq I_1 \cap I_2\). Furthermore, given that \(Z(L) \subseteq I_1 \cap I_2 \subseteq Z(L)\), we conclude that \(I_1 \cap I_2 = Z(L)\).

Next, we examine the structure of \(I_2\). Since \(I_2/Z(L) \cong A(n-5)\), it follows that \(I_2^2 \subseteq Z(L) \cong A(2)\). Consequently, $I_2$  is a nilpotent Lie algebra of nilpotency class at most $2$ and dimension greater than $2.$  We now consider the following cases:
 
\begin{itemize}
    \item[(a)] Suppose  \(c(I_2) = 1\), so that $I_2$ is abelian.Then
\([I_1, I_2] \subseteq I_1 \cap I_2 = Z(L) = Z(I_1) \cong A(2)\).
Since \(\dim I_2 \geq 3\), it follows that $[I_1, I_2] \neq 0 .$ Then $I_2 $ decomposes as
\(I_2 = Z(L) \oplus A\), where \(A \cong A(n-5)\) and  $ [I_1, A] \subseteq Z(L)$ with \(0 \neq [I_1, A] \).
 It follows that

    \[
    L = I_1 + I_2 = I_1 + Z(L) + A = I_1 \rtimes A.
    \]
    
    This corresponds to case $(i).$ 
       \item[(b)] Suppose \(c(I_2) = 2\) and \(I_2^2 \cong A(1)\). Since 
    \[
    \dim I_2 = \dim L - \dim I_1 + \dim(I_1 \cap I_2) = n-5+2 = n-3,
    \]
    by Theorem \ref{Theorem 3.6}, we have a decomposition  $I_2=T_1 \oplus A_1,$
    where \(T_1 \cong H(m)\), \(A_1 \cong A(n-2m-4)\), and \([I_1,I_2] \subseteq I_1 \cap I_2=Z(L)\).
    
    Note that both $[I_1, T_1]$  and $[I_1, A_1]$ lie in the center. Moreover,  $I_1^2=T_1^2=A(1)\subseteq Z(L),$  it follows that $Z(L)=T_1^2\oplus \langle z\rangle$ and   $Z(L)\cap T_1=T_1^2.$ Hence,  $T_1\cap (A_1+\langle z\rangle)=0.$ 
               
     It is straightforward to verify that  $Z(L)\cap A_1=\langle z\rangle.$ Therefore, $\dim A_1\geq 1$ and we can write $A_1=A_2\oplus \langle z\rangle,$\\
where $A_2=A(n-2m-5).$  Hence $I_2=T_1 \oplus A(n-2m-5)\oplus \langle z\rangle.$
   We now divide the analysis into three cases:
    \begin{itemize}
        \item[(b-1)] Let  \([I_1, T_1] = 0\). 
        First, assume that \(A(n-2m-4) = \langle z\rangle\). Then \(n = 2m+5\) and  we get
    \[
    L = I_1 \dot{+} T_1,
    \]
      This corresponds to case $(ii).$
Now, assume that \(A(n-2m-4) \neq A(1)\). Then \(n \neq 2m+5\), and hence
       \[
    L = I_1 + (T_1 \oplus A_2\oplus \langle z\rangle) =( I_1 +A_2+\langle z\rangle ) \dot{+} T_1=(I_1 \rtimes A_2) \dot{+} T_1
    \]
        This corresponds to case $(iii).$
    \item[(b-2)] Let  \([I_1, T_1] =I_1 \cap I_2 =Z(L) \).First, assume that \(A(n-2m-4) = \langle z\rangle \). Then \(n = 2m+5\) and \(I_2 \cong T_1\oplus \langle z\rangle\). Therefore, we get   $ L = I_1 + T_1.$  This corresponds to case $(iv).$\\
Now, assume that \(A(n-2m-4) \neq A(1)\). Then \(n \neq 2m+5\), and hence
    \[
    L = I_1 + (T_1 \oplus A_2 \oplus \langle z\rangle),
    \]
    Since,  $A_2\cap Z(L)=0 $ and   \([I_1, A_1] \neq 0\), we have  $    L =( I_1\rtimes A_2) + T_1.$  It proves case $(v).$
         \item[(b-3)]$[I_1,T_1]\cong A(1).$   First, assume that \(A(n-2m-4) = \langle z\rangle\). Then \(n = 2m+5\) and \(I_2 \cong T_1\oplus \langle z\rangle \). Then $
    L = I_1 + T_1,$ 
    where $[I_1,T_1]=A(1).$ It proves $vi.$\\

Now, assume that \(A(n-2m-4) \neq A(1)\).Then it follows that $n-2m-4\neq 1.$ Consequently, the Lie algebra $L$  can be expressed as $    L = I_1 + (T_1 \oplus A_2 \oplus \langle z\rangle).$
By analyzing all possible relations among $[I_1,T_1], [I_1,A_1], T_1^2$ and the center $Z(L),$ we conclude that the Lie algebra admits a unique structure given by $L = I_1 + (T_1 \oplus A_2 \oplus \langle z\rangle)=(I_1 + T_1) \rtimes A_2. $
It proves $ vii.$
 \item[(c)] Suppose \(c(I_2) = 2\) and \(I_2^2 \cong A(2)\).  
 Since 
    \[
    \dim I_2 = \dim L - \dim I_1 + \dim(I_1 \cap I_2) = n-5+2 = n-3,
    \]
    by \cite[Proposition 2.4]{NP2}, it follows that \(I_2 \cong K\oplus A(t)\),    where \(K\) is a generalized Hiesnberg of rank two. Moreover, since $K^2=I_2^2 \subseteq Z(L),$  and $\dim I_2^2=\dim Z(L),$  we conclude that $ K^2=Z(L)$ and  $Z(L)\cap A(t)=0.$ 
Additionally, as $ [I_1,K] \subseteq Z(L)$  it follows that $K$ is an ideal. Now, if \(A(t) =0\), then $\dim K=n-3$ and $    L = I_1 \dot{+} K,$  which establishes statement $viii.$ 

    On the other hand, if \(A(t) \neq 0\), then $    L = I_1 + (K \oplus A(t)). $  
  In this case, the exact structure of $ L $ can be determined by analyzing the structure of $[I_1,K].$ 
This completes the proof of statement $(ix).$

\end{itemize}

\end{itemize}

   }
 \end{p a}

\begin{p b}
{\em 
By Lemma \ref{Lemma 2.1},  the center factor of $L$ decomposes as 
\[
\frac{L}{Z(L)}= H(1)\oplus A(2).
\]
Assume that
\[
Z(L) = \langle x_{6}, x_{7} \rangle.
\]
Then, the set \(\{ x_{1},\dots,x_{7}\}\) forms a basis for \(L\), and the Lie algebra admits the following structure:
\begin{align*}
L =\langle &x_{1},\dots,x_{7}\mid[x_{1},x_{2}]=x_{3}+\alpha_{1}x_{6}+\alpha_{2}x_{7},~[x_{1},x_{3}]=\alpha_{3}x_{6}+\alpha_{4} x_{7},~[x_{1},x_{4}]=\alpha_{5}\\
&x_{6}+\alpha_{6}x_{7},~[x_{1},x_{5}]=\alpha_{7}x_{6}+\alpha_{8}x_{7},~[x_{2},x_{3}]=\alpha_{9}x_{6}+\alpha_{10}x_{7},~[x_{2},x_{4}]=\alpha_{11}x_{6}+\\
&\alpha_{12}x_{7},~[x_{2},x_{5}]=\alpha_{13}x_{6}+\alpha_{14}x_{7},~[x_{3},x_{4}]=\alpha_{15}x_{6}+\alpha_{16}x_{7},~[x_{3},x_{5}]=\alpha_{17}x_{6}+\\
&\alpha_{18}x_{7},~[x_{4},x_{5}]=\alpha_{19}x_{6}+\alpha_{20}x_{7} \rangle.
\end{align*}
Next, we determine which of the coefficients \(\alpha_i\) vanish. By applying an appropriate change of variables, we may assume  
 \(\alpha_{1} = \alpha_{2} = 0\).
Using the Jacobi identity on the triples \((x_{1},x_{2},x_{4})\) and \((x_{1}, x_{2}, x_{5})\), we deduce that
\[
[x_{3},x_{5}] = 0 \quad \text{and} \quad [x_{3},x_{4}] = 0.
\]
It follows that
\[
\alpha_{15} = \alpha_{16} = \alpha_{17} = \alpha_{18} = 0.
\]
Hence, the Lie algebra \(L\) has the following structure:
\begin{align*}
L =\langle &x_{1},\dots,x_{7}\mid [x_{1},x_{2}]=x_{3},~[x_{1},x_{3}]=\alpha_{1}x_{6}+\alpha_{2}x_{7},~[x_{1},x_{4}]=\alpha_{3}x_{6}+\alpha_{4}x_{7},\\
&[x_{1},x_{5}]=\alpha_{5}x_{6}+\alpha_{6}x_{7},~[x_{2},x_{3}]=\alpha_{7}x_{6}+\alpha_{8}x_{7},~[x_{2},x_{4}]=\alpha_{9}x_{6}+\alpha_{10}x_{7},\\
&[x_{2},x_{5}]=\alpha_{11}x_{6}+\alpha_{12}x_{7},~[x_{4},x_{5}]=\alpha_{13}x_{6}+\alpha_{14}x_{7} \rangle.
\end{align*}

Since \(\dim L^{3} = 2, L^3=[L,L^2]\), and $L^2=\langle x_3,x_6,x_7\rangle$ it follows that the structure constants \(\alpha_{1}\), \(\alpha_{2}\), \(\alpha_{7}\), and \(\alpha_{8}\) cannot be  zero simultaneously.  
Without loss of generality, we may assume
\[
[x_{1},x_{3}] = x_{6} \quad \text{and} \quad [x_{2},x_{3}] = x_{7}.
\]
Under these assumptions, \(L\) has the following structure:
\begin{align*}
L = \langle &x_{1},\dots,x_{7} \mid [x_{1},x_{2}]=x_{3},~[x_{1},x_{3}]=x_{6},~[x_{2},x_{3}]=x_{7},~[x_{1},x_{4}]=\alpha_{1}x_{6}+\alpha_{2}x_{7},\\
&[x_{1},x_{5}]=\alpha_{3}x_{6}+\alpha_{4}x_{7},~[x_{2},x_{4}]=\alpha_{5}x_{6}+\alpha_{6}x_{7},~[x_{2},x_{5}]=\alpha_{7}x_{6}+\alpha_{8}x_{7},\\
&[x_{4},x_{5}]=\alpha_{9}x_{6}+\alpha_{10}x_{7} \rangle.
\end{align*}

To identify which of the \(\alpha_{i}\) are nonzero, we utilize the Schur multiplier of \(L\). 
 Using the method of Stitzinger and Hardy, as described in \cite{H}, we  compute the Schur multiplier of \(L\) and  obtain:
\begin{align*}
&[x_{1},x_{2}]=x_{3}+s_{1},~[x_{1},x_{3}]=x_{6}+s_{2},~[x_{1},x_{4}]=\alpha_{1}x_{6}+\alpha_{2}x_{7}+s_{3},~[x_{1},x_{5}] = \alpha_{3}x_{6}\\
&+\alpha_{4}x_{7}+s_{4},~[x_{1},x_{6}]=s_{5},~[x_{1},x_{7}] = s_{6},~[x_{2},x_{3}]=x_{7}+s_{7},~[x_{2},x_{4}]=\alpha_{5}x_{6}+\alpha_{6}x_{7}\\
&+s_{8},~[x_{2},x_{5}]=\alpha_{7}x_{6}+\alpha_{8}x_{7}+s_{9},~[x_{2},x_{6}]=s_{10},~[x_{2},x_{7}]=s_{11},~[x_{3},x_{4}]=s_{12},\\
&[x_{3},x_{5}]=s_{13},~[x_{3},x_{6}]=s_{14},~[x_{3},x_{7}]=s_{15},~[x_{4},x_{5}]=\alpha_{9}x_{6}+\alpha_{10}x_{7}+s_{16},~[x_{4},x_{6}]\\
&=s_{17},~[x_{4},x_{7}]=s_{18},~[x_{5},x_{6}]=s_{19},~[x_{5},x_{7}]=s_{20},~[x_{6},x_{7}]=s_{21}.
\end{align*}
Without loss of generality, we may assume that \(s_{1} = s_{2} = s_{7} = 0\).
\begin{align*}
&[x_{1},x_{2}]=x_{3},~[x_{1},x_{3}]=x_{6},~[x_{1},x_{4}]=\alpha_{1}x_{6}+\alpha_{2}x_{7}+s_{1},~[x_{1},x_{5}]=\alpha_{3}x_{6}+\alpha_{4}x_{7}+s_{2}\\
&[x_{1},x_{6}]=s_{3},~[x_{1},x_{7}]=s_{4},~[x_{2},x_{3}]=x_{7},~[x_{2},x_{4}]=\alpha_{5}x_{6}+\alpha_{6}x_{7}+s_{5},~[x_{2},x_{5}]=\alpha_{7}x_{6}\\
&+\alpha_{8}x_{7}+s_{6},~[x_{2},x_{6}]=s_{7},~[x_{2},x_{7}]=s_{8},~[x_{3},x_{4}]=s_{9},~[x_{3},x_{5}]=s_{10},~[x_{3},x_{6}]=s_{11}\\
&[x_{3},x_{7}]=s_{12},~[x_{4},x_{5}]=\alpha_{9}x_{6}+\alpha_{10}x_{7}+s_{13},~[x_{4},x_{6}]=s_{14},~[x_{4},x_{7}]=s_{15},~[x_{5},x_{6}]=\\
&s_{16},~[x_{5},x_{7}]=s_{17},~[x_{6},x_{7}]=s_{18}.
\end{align*}
Thus, we have \(M(L) = \langle s_{1},\dots,s_{18} \rangle\).
Applying the Jacobi identity once more, we derive the following relations:
\begin{align*}
&s_{4} = s_{7},\quad \text{and} \quad s_{11} = s_{12} = s_{14} = s_{15} = s_{16} = s_{17} = s_{18} = 0,\\
&\alpha_{5}s_{3}+\alpha_{6}s_{4}-\alpha_{1}s_{7}-\alpha_{2}s_{8}-s_{9}=0,\quad\alpha_{7}s_{3} + \alpha_{8}s_{4} - \alpha_{3}s_{7} - \alpha_{4}s_{8} - s_{10} = 0\\
&\alpha_{9}s_{3} + \alpha_{10}s_{4} = 0,\quad \alpha_{9}s_{7} + \alpha_{10}s_{8} = 0.
\end{align*}
Therefore, \(M(L) = \langle s_{1}, s_{2}, s_{3}, s_{4}, s_{5}, s_{6}, s_{8}, s_{13} \rangle\).
Upon further examination, the above relations give rise to several possible cases.\\
\textbf{Case (1):} If \(s_{3},~s_{4}\), and \(s_{8}\) are nonzero and linearly independent, then the Schur multiplier has dimension eight. Consequently, \(\alpha_{9}=\alpha_{10}=0\).
Therefore, the Lie algebra \(L\) admits the following presentation:
\begin{align*}
L = \langle &x_{1},\dots,x_{7} \mid [x_{1},x_{2}]=x_{3},~[x_{1},x_{3}]=x_{6},~[x_{2},x_{3}]=x_{7},~[x_{1},x_{4}]=\alpha_{1}x_{6}+\alpha_{2}x_{7},\\
&[x_{1},x_{5}]=\alpha_{3}x_{6}+\alpha_{4}x_{7},~[x_{2},x_{4}]=\alpha_{5}x_{6}+\alpha_{6}x_{7},~[x_{2},x_{5}]=\alpha_{7}x_{6}+\alpha_{8}x_{7} \rangle.
\end{align*}
\textbf{Subcase 1.1}
Assume that \(\alpha_{1}=\alpha_{6}\) and \(\alpha_{3}=\alpha_{8}\).
We claim that under these conditions, \(L\) is isomorphic to \(L_{2}\).
To prove this, we proceed step by step.
First, replace \(x_{4}\) with \(x_{4}-\alpha_{1}x_{3}\) and \(x_{5}\) with \(x_{5}-\alpha_{3}x_{3}\).
This transformation simplifies the structure of \(L\) to:
\begin{align*}
L = \langle x_{1},\dots,x_{7} \mid &[x_{1},x_{2}]=x_{3},~[x_{1},x_{3}]=x_{6},~[x_{2},x_{3}]=x_{7},~[x_{1},x_{4}]=\alpha_{2}x_{7}\\
&[x_{1},x_{5}]=\alpha_{4}x_{7},~[x_{2},x_{4}]=\alpha_{5}x_{6},~[x_{2},x_{5}]=\alpha_{7}x_{6} \rangle.
\end{align*}
Next, replace \(x_{4}\) with \(\alpha_{2}^{-1}x_{4}\) and \(x_{5}\) with \(\alpha_{4}^{-1}x_{5}\).
After this change of basis, we obtain:
\begin{align*}
L = \langle x_{1},\dots,x_{7} \mid&[x_{1},x_{2}]=x_{3},~[x_{1},x_{3}]=x_{6},~[x_{2},x_{3}]=x_{7},~[x_{1},x_{4}]=x_{7}\\&[x_{1},x_{5}]=x_{7},~[x_{2},x_{4}]=\alpha_{5}x_{6},~[x_{2},x_{5}]=\alpha_{7}x_{6} \rangle.
\end{align*}
Then replace \(x_{5}\) with \(x_{5}-x_{4}\). As a result, we have:
\begin{align*}
L = \langle x_{1},\dots,x_{7} \mid &[x_{1},x_{2}]=x_{3},~[x_{1},x_{3}]=x_{6},~[x_{2},x_{3}]=x_{7},~[x_{1},x_{4}]=x_{7}\\
&[x_{2},x_{4}]=\alpha_{5}x_{6},~[x_{2},x_{5}]=\alpha_{7}x_{6} \rangle.
\end{align*}

Continuing the process, rename \(x_{5}\) as \(\alpha_{7}^{-1}x_{5}\) to get:
\begin{align*}
L =\langle x_{1},\dots,x_{7} \mid &[x_{1},x_{2}]=x_{3},~[x_{1},x_{3}]=x_{6},~[x_{2},x_{3}]=x_{7},~[x_{1},x_{4}]=x_{7}\\
&[x_{2},x_{4}]=\alpha_{5}x_{6},~[x_{2},x_{5}]=x_{6} \rangle.
\end{align*}
Finally, replacing \(x_{4}\) with \(x_{4}-\alpha_{5}^{-1}x_{5}\), we arrive at the presentation:
\begin{align*}
L = \langle x_{1},\dots,x_{7} \mid &[x_{1},x_{2}]=x_{3},~[x_{1},x_{3}]=x_{6},~[x_{2},x_{3}]=x_{7},~[x_{1},x_{4}]=x_{7},~[x_{2},x_{5}]=x_{6} \rangle.
\end{align*}

It is clear that by swapping \(x_{4}\) and \(x_{5}\), the Lie algebra becomes isomorphic to \(L_{2}\).

\item[\textbf{Subcase 1.2:}]  If we assume that \(\alpha_{2} = \alpha_{5}\), \(\alpha_{3} = \alpha_{8}\), and \(\alpha_{4} = \alpha_{7}\), then \(L\) is isomorphic to \(L_{1}\).

To prove this claim, we examine the structure of $ L $ under these assumptions:
\begin{align*}
L = \langle &x_{1},\dots,x_{7} \mid [x_{1},x_{2}]=x_{3},~[x_{1},x_{3}]=x_{6},~[x_{2},x_{3}]=x_{7},~[x_{1},x_{4}]=\alpha_{1}x_{6}+\alpha_{2}x_{7}\\
&[x_{1},x_{5}]=\alpha_{3}x_{6}+\alpha_{4}x_{7},~[x_{2},x_{4}]=\alpha_{2}x_{6}+\alpha_{6}x_{7},~[x_{2},x_{5}]=\alpha_{4}x_{6}+\alpha_{3}x_{7} \rangle.
\end{align*}
First, replace \(x_{5}\) with \(x_{5}-\alpha_{3}x_{3}\), so that:
\begin{align*}
L=\langle x_{1},\dots,x_{7} \mid &[x_{1},x_{2}]=x_{3},~[x_{1},x_{3}]=x_{6},~[x_{2},x_{3}]=x_{7},~[x_{1},x_{4}]=\alpha_{1}x_{6}+\alpha_{2}x_{7}\\
&[x_{1},x_{5}]=\alpha_{4}x_{7},~[x_{2},x_{4}]=\alpha_{2}x_{6}+\alpha_{6}x_{7},~[x_{2},x_{5}]=\alpha_{4}x_{6} \rangle.
\end{align*}
Next, replace \(x_{4}\) with \(\alpha_{2}^{-1}x_{4}\) and \(x_{5}\) with \(\alpha_{4}^{-1}x_{5}\). Thus, we have:
\begin{align*}
L = \langle x_{1},\dots,x_{7} \mid &[x_{1},x_{2}]=x_{3},~[x_{1},x_{3}]=x_{6},~[x_{2},x_{3}]=x_{7},~[x_{1},x_{4}]=\alpha_{1}x_{6}+x_{7}\\
&[x_{1},x_{5}]=x_{7},~[x_{2},x_{4}]=x_{6}+\alpha_{6}x_{7},~[x_{2},x_{5}]=x_{6} \rangle.
\end{align*}

Then, replace \(x_{4}\) with \(x_{4}-x_{5}\), yielding:
\begin{align*}
L = \langle x_{1},\dots,x_{7} \mid &[x_{1},x_{2}]=x_{3},~[x_{1},x_{3}]=x_{6},~[x_{2},x_{3}]=x_{7},~[x_{1},x_{4}]=\alpha_{1}x_{6}\\&[x_{1},x_{5}]=x_{7},~[x_{2},x_{4}]=\alpha_{6}x_{7},~[x_{2},x_{5}]=x_{6} \rangle.
\end{align*}
Next, replace \(x_{4}\) with \(\alpha_{1}^{-1}x_{4}\). Therefore, we obtain:
\begin{align*}
L = \langle x_{1},\dots,x_{7} \mid &[x_{1},x_{2}]=x_{3},~[x_{1},x_{3}]=x_{6},~[x_{2},x_{3}]=x_{7},~[x_{1},x_{4}]=x_{6}\\&[x_{1},x_{5}]=x_{7},~[x_{2},x_{4}]=\alpha_{6}x_{7},~[x_{2},x_{5}]=x_{6} \rangle.
\end{align*}

Now, replace \(x_{4}\) with \(x_{4}-\alpha_{6}x_{3}\), resulting in:
\begin{align*}
L = \langle x_{1},\dots,x_{7} \mid& [x_{1},x_{2}]=x_{3},~[x_{1},x_{3}]=x_{6},~[x_{2},x_{3}]=x_{7},~[x_{1},x_{4}]=(1-\alpha_{6})x_{6}\\&[x_{1},x_{5}]=x_{7},~[x_{2},x_{5}]=x_{6} \rangle.
\end{align*}

Afterwards, replace \(x_{4}\) with \((1-\alpha_{6})^{-1}x_{4}\) and replace \(x_{2}\) with \(x_{2}-x_{1}\), yielding:
\begin{align*}
L=\langle x_{1},\dots,x_{7} \mid& [x_{1},x_{2}]=x_{3},~[x_{1},x_{3}]=x_{6},~[x_{2},x_{3}]=x_{7}-x_{6},~[x_{1},x_{4}]=x_{6}\\&[x_{1},x_{5}]=x_{7},~[x_{2},x_{5}]=x_{6}-x_{7} \rangle.
\end{align*}

Finally, by replacing \(x_{5}\) with \(x_{5}+x_{3}-2x_{4}\) and renaming \(x_{7}-x_{6}\) as \(x_{7}\), the resulting Lie algebra is isomorphic to \(L_{1}\).\\
\textbf{Subcase 1.3:} Suppose that \(\alpha_{1} = \alpha_{6}\), \(\alpha_{2} = \alpha_{5}\), and \(\alpha_{4} = \alpha_{7}\). Then the Lie algebra $L$ exhibits the following structure:
\begin{align*}
L=\langle &x_{1},\dots,x_{7} \mid [x_{1},x_{2}]=x_{3},~[x_{1},x_{3}]=x_{6},~[x_{2},x_{3}]=x_{7},~[x_{1},x_{4}]= \alpha_{1}x_{6}+\alpha_{2}x_{7}\\&[x_{1},x_{5}]=\alpha_{3}x_{6}+\alpha_{4}x_{7},~
[x_{2},x_{4}]=\alpha_{2}x_{6}+\alpha_{1}x_{7},~[x_{2},x_{5}]=\alpha_{4}x_{6}+\alpha_{8}x_{7}\rangle.
\end{align*}
By replacing \(x_{4}\) with \(x_{4} - \alpha_{1}x_{3}\) and \(x_{5}\) by \(\alpha_{4}^{-1}x_{5}\), the Lie algebra takes the following structure:
\begin{align*}
L = \langle &x_{1},\dots,x_{7} \mid [x_{1},x_{2}]=x_{3},~[x_{1},x_{3}]=x_{6},~[x_{2},x_{3}]=x_{7},~[x_{1},x_{4}] = \alpha_{2}x_{7}\\&[x_{1},x_{5}]=\alpha_{3}x_{6}+x_{7},~[x_{2},x_{4}] = \alpha_{2}x_{6},~[x_{2},x_{5}] = x_{6} + \alpha_{8}x_{7} \rangle.
\end{align*}
Now, first replace \(x_{4}\) with \(\alpha_{2}^{-1}x_{4}\), and then replace \(x_{5}\) with \(x_{5} - x_{4}\). We obtain:
\begin{align*}
L = \langle x_{1},\dots,x_{7} \mid &[x_{1},x_{2}]=x_{3},~[x_{1},x_{3}]=x_{6},~[x_{2},x_{3}]=x_{7},~[x_{1},x_{4}] = x_{7}\\&[x_{1},x_{5}]=\alpha_{3}x_{6},~[x_{2},x_{4}]=x_{6},~[x_{2},x_{5}]=\alpha_{8}x_{7}\rangle.
\end{align*}
Next, replace \(x_{5}\) with \(x_{5} - \alpha_{3}x_{3}\). Therefore, we have:
\begin{align*}
L = \langle x_{1},\dots,x_{7} \mid &[x_{1},x_{2}] = x_{3},~[x_{1},x_{3}] = x_{6},~[x_{2},x_{3}] = x_{7},~[x_{1},x_{4}] = x_{7}\\&[x_{2},x_{4}]=x_{6},~[x_{2},x_{5}] = \alpha_{8}x_{7} \rangle.
\end{align*}
Now, in the first step, replace \(x_{5}\) with \(\alpha_{8}^{-1}x_{5}\).
In the second step, replace \(x_{1}\) with \(x_{1} - x_{2}\).
In the third step, replace \(x_{4}\) with \(x_{4} + x_{3}\).
Finally, rename all \(x_{6} - x_{7}\) as \(x_{6}\).
Thus we obtain:
\[
L = \langle x_{1},\dots,x_{7} \mid [x_{1},x_{2}] = x_{3},~[x_{1},x_{3}] = x_{6},~[x_{2},x_{3}] = x_{7},~[x_{2},x_{4}] = x_{6},~[x_{2},x_{5}] = x_{7} \rangle.
\]

Finally, applying the change of variables, we arrive at:
\[
x_{1}' = x_{2},\quad x_{2}' = x_{1},\quad x_{3}' = -x_{3},\quad x_{4}' = -x_{5},\quad x_{5}' = -x_{4},\quad x_{6}' = -x_{7},\quad x_{7}' = -x_{6},
\]
the resulting Lie algebra is isomorphic to \(L_{1}\).\\
\textbf{Case (2):} Assume that \(\alpha_{10} = 0\) and \(\alpha_{9} \neq 0\), then $s_3=s_4=0$, \(\dim M(L) = 6\) and \(L\) has the following structure:
\begin{align*}
L = \langle x_{1},\dots,x_{7}\mid&[x_{1},x_{2}]=x_{3},~[x_{1},x_{3}]=x_{6},~[x_{2},x_{3}]=x_{7},~[x_{1},x_{4}]=\alpha_{1}x_{6} + \alpha_{2}x_{7}\\&[x_{1},x_{5}]=\alpha_{3}x_{6}+\alpha_{4}x_{7},~[x_{2},x_{4}]=\alpha_{5}x_{6}+\alpha_{6}x_{7}\\&[x_{2},x_{5}] = \alpha_{7}x_{6}+\alpha_{8}x_{7},~[x_{4},x_{5}]=\alpha_{9}x_{6} \rangle.
\end{align*}
\textbf{Subcase 2.1:} If, in this case, we also assume that \(\alpha_{1} = \alpha_{6}\), \(\alpha_{2} = \alpha_{4}\), \(\alpha_{3} = \alpha_{8}\), and \(\alpha_{5} = \alpha_{7}\), then we claim that \(L\) is isomorphic to \(L_{4}\).
To prove this claim, we proceed as follows:
First, replace \(x_{4}\) with \(x_{4} - \alpha_{1}x_{3}\) and \(x_{5}\) with \(x_{5} - \alpha_{3}x_{3}\).
Thus, the Lie algebra \(L\) becomes:
\begin{align*}
L = \langle x_{1},\dots,x_{7}\mid&[x_{1},x_{2}]=x_{3},~[x_{1},x_{3}]=x_{6},~[x_{2},x_{3}]=x_{7},~[x_{1},x_{4}] = \alpha_{2}x_{7}\\&[x_{1},x_{5}]=\alpha_{2}x_{7},~[x_{2},x_{4}]=\alpha_{5}x_{6},~[x_{2},x_{5}]=\alpha_{5}x_{6},~[x_{4},x_{5}]=\alpha_{9}x_{6} \rangle.
\end{align*}
Next, replace \(x_{5}\) with \(x_{5} - x_{4}\).
Then, rename \(\alpha_{2}^{-1}x_{4}\) as \(x_{4}\) and \((\alpha_{2}\alpha_{9})^{-1}x_{5}\) as \(x_{5}\).
We obtain:
\begin{align*}
L = \langle x_{1},\dots,x_{7}\mid&[x_{1},x_{2}]=x_{3},~[x_{1},x_{3}]=x_{6},~[x_{2},x_{3}]=x_{7},~[x_{1},x_{4}] = x_{7}\\&[x_{2},x_{4}]=\alpha_{5}x_{6},~[x_{4},x_{5}]=x_{6}\rangle.
\end{align*}

Finally, replace \(x_{2}\) with \(x_{2} - x_{4}\), \(x_{3}\) with \(x_{3} + x_{7}\), and \(x_{4}\) with \(x_{4} - x_{5}\).
Thus, \(L\) becomes:
\[
L = \langle x_{1}, \dots, x_{7} \mid [x_{1},x_{2}] = x_{3},~[x_{1},x_{3}] = x_{6},~[x_{2},x_{3}] = x_{7},~[x_{1},x_{4}] = x_{7},~[x_{4},x_{5}] = x_{6} \rangle.
\]

Now, it is obvious that by renaming 
\[
x_{1}' = x_{2}, \quad x_{2}' = x_{1}, \quad x_{3}' = -x_{3}, \quad x_{6}' = -x_{7}, \quad x_{7}' = -x_{6},
\]
the Lie algebra \(L\) is isomorphic to \(L_{4}\).

\textbf{Case (iii):} If \(\alpha_{9} = 0\) and \(\alpha_{10} \neq 0\), then $s_4=s_8=0$, \(\dim M(L) = 6\) and \(L\) has the following structure:
\begin{align*}
L = \langle x_{1},\dots,x_{7}\mid&[x_{1},x_{2}]=x_{3},~[x_{1},x_{3}]=x_{6},~[x_{2},x_{3}]=x_{7},~[x_{1},x_{4}]=\alpha_{1}x_{6} + \alpha_{2}x_{7}\\&[x_{1},x_{5}]=\alpha_{3}x_{6}+\alpha_{4}x_{7},~[x_{2},x_{4}]=\alpha_{5}x_{6}+\alpha_{6}x_{7}\\&[x_{2},x_{5}] = \alpha_{7}x_{6}+\alpha_{8}x_{7},~[x_{4},x_{5}]=\alpha_{10}x_{7}\rangle.
\end{align*}
\textbf{Subcase 3.1:} If, in this case, we also assume that \(\alpha_{1} = \alpha_{6}\), \(\alpha_{2} = \alpha_{5}\), and \(\alpha_{4} = \alpha_{7}\), then the Lie algebra \(L\) will be isomorphic to \(L_{4}\).

To prove this claim, we proceed step by step. 

First, replace \(x_{4}\) with \(x_{4} - \alpha_{1}x_{3}\), and \(x_{5}\) with \(x_{5} - \alpha_{4}\alpha_{2}^{-1}x_{4}\). We obtain:
\begin{align*}
L = \langle &x_{1}, \dots, x_{7} \mid [x_{1},x_{2}]=x_{3},~[x_{1},x_{3}]=x_{6},~[x_{2},x_{3}]=x_{7},~[x_{1},x_{4}]=\alpha_{2}x_{7}\\&[x_{1},x_{5}]=\alpha_{3}x_{6},~[x_{2},x_{4}]=\alpha_{2}x_{6},~[x_{2},x_{5}]=\alpha_{8}x_{7},~[x_{4},x_{5}]=\alpha_{10}x_{7} \rangle.
\end{align*}
Next, replace \(x_{5}\) with \(x_{5} - \alpha_{3}x_{3}\). Then the Lie algebra becomes:
\begin{align*}
L = \langle x_{1},\dots,x_{7}\mid&[x_{1},x_{2}]=x_{3},~[x_{1},x_{3}]=x_{6},~[x_{2},x_{3}]= x_{7},~[x_{1},x_{4}] = \alpha_{2}x_{7}\\&[x_{2},x_{4}]=\alpha_{2}x_{6},~[x_{2},x_{5}]=(\alpha_{8}-\alpha_{3})x_{7},~[x_{4},x_{5}]=\alpha_{10}x_{7} \rangle.
\end{align*}

Now, rename \(\alpha_{2}^{-1}x_{4}\) as \(x_{4}\) and \((\alpha_{8} - \alpha_{3})^{-1}x_{5}\) as \(x_{5}\), thus we obtain:
\begin{align*}
L = \langle x_{1}, \dots, x_{7} \mid &[x_{1},x_{2}] = x_{3},~[x_{1},x_{3}] = x_{6},~[x_{2},x_{3}] = x_{7},~[x_{1},x_{4}] = x_{7}\\&[x_{2},x_{4}]=x_{6},~[x_{2},x_{5}] = x_{7},~[x_{4},x_{5}] = \alpha_{10}x_{7} \rangle.
\end{align*}

In the next step, first replace \(x_{2}\) with \(x_{2} - \alpha_{10}^{-1}x_{4}\), and then rename \(x_{3} + \alpha_{10}x_{7}\) as \(x_{3}\). Therefore, \(L\) becomes:
\begin{align*}
L = \langle x_{1}, \dots, x_{7} \mid &[x_{1},x_{2}] = x_{3},~[x_{1},x_{3}] = x_{6},~[x_{2},x_{3}] = x_{7},~[x_{1},x_{4}] = x_{7}\\&[x_{2},x_{4}] = x_{6},~[x_{4},x_{5}] = \alpha_{10}x_{7} \rangle.
\end{align*}

Finally, rename \(x_{5}\) as \(\alpha_{10}^{-1}x_{5}\) and replace \(x_{1}\) with \(x_{1} + x_{5}\). Thus, the resulting Lie algebra is isomorphic to \(L_{4}\).\\
\textbf{Subcase 3.2:} Assume that \(\alpha_{1} = \alpha_{6}\), \(\alpha_{2} = \alpha_{5}\), \(\alpha_{3} = \alpha_{8}\), and \(\alpha_{4} = \alpha_{7}\). Then \(L\) has the following structure:
\begin{align*}
L =\langle &x_{1}, \dots, x_{7} \mid [x_{1},x_{2}] = x_{3},~[x_{1},x_{3}] = x_{6},~[x_{2},x_{3}] = x_{7},~[x_{1},x_{4}] = \alpha_{1}x_{6} + \alpha_{2}x_{7}\\&[x_{1},x_{5}] = \alpha_{3}x_{6} + \alpha_{4}x_{7},~[x_{2},x_{4}] = \alpha_{2}x_{6} + \alpha_{1}x_{7},~[x_{2},x_{5}] = \alpha_{4}x_{6} + \alpha_{3}x_{7},\\
&[x_{4},x_{5}] = \alpha_{10}x_{7} \rangle.
\end{align*}
First, replace \(x_{4}\) with \(x_{4} - \alpha_{1}x_{3}\) and \(x_{5}\) with \(x_{5} - \alpha_{3}x_{3}\). Then we obtain:
\begin{align*}
L = \langle &x_{1},\dots,x_{7}\mid[x_{1},x_{2}]=x_{3},~[x_{1},x_{3}]=x_{6},~[x_{2},x_{3}]=x_{7},~[x_{1},x_{4}]=\alpha_{2}x_{7}\\&[x_{1},x_{5}]=\alpha_{4}x_{7},~[x_{2},x_{4}]=\alpha_{2}x_{6},~[x_{2},x_{5}]=\alpha_{4}x_{6},~[x_{4},x_{5}]=\alpha_{10}x_{7} \rangle.
\end{align*}

Next, rename \(\alpha_{2}^{-1}x_{4}\) as \(x_{4}\) and \(\alpha_{4}^{-1}x_{5}\) as \(x_{5}\), thus:
\begin{align*}
L =\langle x_{1},x_{2},\dots,x_{7}\mid&[x_{1},x_{2}]=x_{3},~[x_{1},x_{3}] = x_{6},~[x_{2},x_{3}] = x_{7},~[x_{1},x_{4}]=x_{7}\\&[x_{1},x_{5}]=x_{7},~[x_{2},x_{4}]=x_{6},~[x_{2},x_{5}]=x_{6},~[x_{4},x_{5}]=\alpha_{10}x_{7}\rangle.
\end{align*}

Then replace \(x_{4}\) with \(x_{4} - x_{5}\) and \(x_{5}\) with \(x_{5} - x_{4}\). We have:
\begin{align*}
L = \langle x_{1},x_{2},\dots,x_{7}\mid[x_{1},x_{2}]=x_{3},~[x_{1},x_{3}]=x_{6},~[x_{2},x_{3}]=x_{7},~[x_{4},x_{5}] = \alpha_{10}x_{7} \rangle.
\end{align*}

Finally, by replacing \(x_{5}\) with \(\alpha_{10}^{-1}x_{5}\), the resulting Lie algebra is isomorphic to \(L_{3}\).}

\end{p b}

{\footnotesize {\bf Saboura Yousefi}\; \\ {Department of
Mathematics}, {University of Birjand, P.O.Box  97175/615,} {Birjand, Iran.}\\
{\tt Email: saboura.yousefi@birjand.ac.ir}\\

{\footnotesize {\bf Azam Kaheni}\; \\ {Department of
Mathematics}, {University of Birjand, P.O.Box  97175/615,} {Birjand, Iran.}\\
{\tt Email: azamkaheni@birjand.ac.ir}\\

{\footnotesize {\bf Farangis Johari}\; \\ {Center of Mathematics, Computer Science and Cognition}, {Federal University of ABC, Av.\ dos Estados, 5001 -- Bang\'u ,Santo Andr\'e,09210-580,} {S\~ao Paulo,Brazil.}\\
{\tt Email: farangis.johari@ufabc.edu.br}\\

\end{document}